\documentclass[reqno]{amsart}
\usepackage{amsfonts}
\usepackage{amssymb}
\usepackage{amsthm}
\usepackage{upref}
\usepackage{enumerate}
\makeatletter
\def\LaTeX{\leavevmode L\raise.42ex
    \hbox{\kern-.3em\size{\sf@size}{0pt}\selectfont A}\kern-.15em\TeX}
 \makeatother
 
\numberwithin{equation}{section}

\newtheorem{lemma}{Lemma}[section]
\newtheorem{theorem}[lemma]{Theorem} 

\newtheorem{proposition}[lemma]{Proposition}

\theoremstyle{definition}

\newtheorem{example}[lemma]{Example}

\newtheorem{problem}[lemma]{Problem}

 \newcommand{\supp}{\operatorname{supp}}
  \newcommand{\e}{\eqref}
\newcommand{\ri}{\rightarrow}
\newcommand{\q}{\quad}
\newcommand{\ii}{\infty}
\newcommand{\h}{\hbar}
\newcommand{\la}{\langle}
\newcommand{\ra}{\rangle}
\renewcommand{\d}{\delta}

\newenvironment{pf}{\begin{proof}}{\end{proof}}

\def\qqq{\mathrel{\subset\mkern-15mu\lower.38ex\hbox{${\scriptscriptstyle\rightarrow}$}}}

\let\Bbb\mathbb

\begin{document}
\title[A balance of  energies]
{A generalized virial theorem and the balance of kinetic and potential energies in the semiclassical limit}
\author{ D. R. Yafaev}
\address{ IRMAR, Universit\'{e} de Rennes I\\ Campus de
  Beaulieu, 35042 Rennes Cedex, FRANCE}
\email{yafaev@univ-rennes1.fr}
\thanks{Partially supported by the project NONAa, ANR-08-BLANC-0228}
 
\dedicatory{To the memory of Pierre Duclos}

\begin{abstract}
 We obtain  two-sided bounds on kinetic and potential energies  of  a  bound  state of a quantum particle in the semiclassical  limit, as the Planck constant $\hbar\ri 0$. 
 Proofs of these results rely on the generalized virial theorem obtained in the paper as well as on a decay of eigenfunctions in the classically forbidden region.
\end{abstract}
\maketitle

\thispagestyle{empty}

\section{Introduction}
{\bf 1.1. }
Let us consider an eigenfunction $\psi_{\h}(x ) $ of the Schr\"odinger operator
\begin{equation}
H_{\h}=-\h^2 \Delta+v(x),\q v(x)=\overline{v(x)},
\label{eq:I1}\end{equation}
in the space $L_{2}({\Bbb R}^d)$, that is 
\begin{equation}
 -\h^2 \Delta \psi_{\h}(x ) +v(x) \psi_{\h}(x ) =\lambda_{\h}\psi_{\h}(x ), \q \psi_{\h}\in L_{2}({\Bbb R}^d).
\label{eq:I2}\end{equation}
We  always suppose that the function  $v(x)$ is semibounded from below, and hence we can set
\[
\min_{x\in {\Bbb R}^d} v(x)=0.
\]
We typically assume that the potential   $v(x)$ contains   wells and study $\psi_{\h}(x ) $ for $\lambda_{\h}$  close to some non-critical value $\lambda_{0} > 0$ (that is $\nabla v(x) \neq 0$ for all $x$ such that $v(x)=\lambda_{0}$). In particular, $\lambda_{\h}$ is separated from   bottoms of   potential wells.
     The eigenfunctions $\psi_{\h}(x )$ are supposed to be real and normalized, that is 
\[
  \int _{{\Bbb R}^d}   \psi_{\h}^2(x )  dx =1.
\]

 Our goal is to study the behavior of the kinetic 
 \[
K(\psi_{\h})= \h^2  \int _{{\Bbb R}^d}   |  \nabla \psi_{\h}(x )|^2 dx 
\]
and potential 
\[
U(\psi_{\h})=   \int _{{\Bbb R}^d} v(x)      \psi_{\h}^2 (x ) dx 
\]
energies as $\h\ri 0$ if the total energy  
\begin{equation}
\lambda_{\h}=K(\psi_{\h}) + U(\psi_{\h}) 
\label{eq:Tot}\end{equation}
 is close to   $\lambda_{0}$.
 
 \medskip
 
 {\bf 1.2.} 
  To be more precise, we discuss the following  
 
 \begin{problem}\label{Q1}
 Is it true that 
 \begin{eqnarray}
  K(\psi_{\h})\geq c>0
 \label{eq:maink}\end{eqnarray}
for all $\lambda_{\h}$ in a neighborhood of a non-critical point $\lambda_{0} > 0$ uniformly in $\h$?
\end{problem}

In view of equation \e{eq:Tot} this problem  can be equivalently reformulated in terms of the potential energy $U(\psi_{\h})$ as the inequality
 \begin{eqnarray}
  U(\psi_{\h})\leq \lambda_{\h}-c.
 \label{eq:vx}\end{eqnarray}
  
Note that, for small $\h$, the eigenfunctions $\psi_{\h}$ are essentially localized (see \cite{ShHe, Simon, Hel})  in the classically allowed region where $v(x)\leq\lambda_{\h}+\varepsilon$ (for an arbitrary $\varepsilon>0$). 
Therefore, inequality \e{eq:vx}, roughly speaking, means that the eigenfunctions $\psi_{\h}$ are not too strongly localized in a neighborhood of the set $v(x)=\lambda_{\h}$.

Let us discuss Problem~\ref{Q1} in a heuristic way. The first level of the discussion is very superficial. Actually, as $\h\ri 0$, one might expect that the term $-\h^2 \Delta$ disappears so that in the limit we obtain the operator of multiplication by the function $v(x)$.  This operator has a continuous spectrum and its ``eigenfunctions" are Dirac functions of the variable $v(x)-\lambda$. Such functions ``live" in a neighborhood of the set $v(x)=\lambda$ which might eventually interfere with the positive answer to Problem~\ref{Q1}. 

The second level is, on the contrary, quite deep and stipulates that, for small $\h$, the behavior of a quantum particle with Hamiltonian \e{eq:I1} is close to the behavior of the corresponding classical particle, and hence the classical equations of motion can be used. In this context we mention book \cite{FedMas} relying on the method of Maslov canonical operator and papers  \cite{ShRo, ShMaRo} relying on  methods  of microlocal analysis. 

We avoid this deep level using only the virial theorem and the fact that a quantum particle of energy $\lambda_{\h}$ should be essentially localized as $\h \to 0$ in the  classically allowed region.

 \medskip
 
 {\bf 1.3.}
Problem~\ref{Q1} arose by the proof of the limiting absorption principle for the   Hamiltonian $K$ of a quantum particle moving in a magnetic field of an infinite straight current (see  \cite{Y}). This problem reduces to a study of   eigenvalues $\lambda_{\h}$ close to a point $\lambda_{0} >0$ for a   potential 
 \[
v(x)=v_{0}\ln^2 |x|, \q v_{0}>0,\q x\in{\Bbb R}^2. 
\]
 The limiting absorption principle for the operator $K$  requires the estimate
 \[
d \lambda_{\h}/ d\h\geq c \h^{-1},\q c>0.  
\]
It looks somewhat exotic but in view of the formula (see, e.g., \cite{Kato})
\[
d \lambda_{\h}/ d\h= 2\h 
\int _{{\Bbb R}^d} |\nabla  \psi_{\h}(x )|^2 dx, 
\]
it is equivalent to estimate \e{eq:maink}.

 \medskip
 
 {\bf 1.4.}
 We also discuss   a problem dual to Problem~\ref{Q1}.  

 \begin{problem}\label{Q2}
  Is it true that 
\[  
  K(\psi_{\h})\leq \lambda_{\h}-c,\q c>0,
\]
for all $\lambda$ in a  neighborhood of a non-critical point $\lambda_{0}>0$ uniformly in  $\h$?
\end{problem}

This fact is equivalent to the estimate $  U(\psi_{\h})\geq c>0$ which   means that, in the semiclassical limit,  the eigenfunctions $\psi_{\h}$ are not too strongly localized   at the bottom of the potential well.

\section{A generalized virial theorem}

 {\bf 2.1.}
 Below the operator $H_{\h}$ is always defined by formula \e{eq:I1}.   The following result generalizes the classical  virial theorem. 
 
  \begin{theorem}\label{GVT}
   Suppose that  $v\in C^1 ({\Bbb R}^d)$. Let $a =\bar{a} \in C^4 ({\Bbb R}^d)$, and let its four derivatives be bounded. Then eigenfunctions $\psi_{\h}$ of the operator $H_{\h}$  satisfy an identity
   \begin{equation}
\int _{{\Bbb R}^d} \Bigl( 4 \h^2\sum_{j,k=1}^d    a_{jk} \partial_{k}\psi_{\h} \partial_{j}\psi_{\h} -
 \h^2  (\Delta^2 a) \psi_{\h}^2 - 
2  \la\nabla a, \nabla v \ra \psi_{\h}^2   \Bigl)dx=0, 
\label{eq:V1M}\end{equation}
where $   a_{jk} =\partial^2 a/\partial x_{j} \partial x_{k}$.
       \end{theorem}
   
   \begin{pf}
   Let  
\[
A=   \sum_{j=1}^d \big(a_{j}(x)  \partial_{j} + \partial_{j} a_{j}(x)\big),\q a_{j} =\partial  a/\partial x_{j},
\]
be a   general self-adjoint  first order differential operator.
Then the commutators
\[
[-\Delta, A]= - 4 \sum_{j,k=1}^d \partial_{j}a_{jk}\partial_{k}-  \Delta^2 a 
\]
and 
\[
[v, A]= -2  \la\nabla a, \nabla v\ra.
\]
It remains to use that, for a function $\psi_{\h}$ satisfying equation \e{eq:I2}, the
identity
\[
([H_{\h}  ,A]\psi_{\h}, \psi_{\h})=0
\]
holds.
       \end{pf}
       
       Note that since eigenfunctions corresponding to isolated eigenvalues decay exponentially, identity \e{eq:V1M} requires practically no assumptions on the behavior of the function $a(x)$ as $|x|\to\infty$.  However we consider only bounded functions $a(x)$.

    If  $a(x)=x^2 $, then identity \e{eq:V1M} reduces to the usual (see, e.g., \cite{RS}) form 
    \begin{equation}
2\h^2 \int _{{\Bbb R}^d}|\nabla \psi_{\h}(x )|^2 dx=\int _{{\Bbb R}^d}
r v_{r}(x)  \psi_{\h}^2 (x ) dx, \q r=|x|,
\label{eq:V1}\end{equation}
of the virial theorem.
      Combining equations  \e{eq:Tot} 
 and \e{eq:V1} we see that
\begin{equation}
  \int _{{\Bbb R}^d} \big( 2^{-1}r v_{r}(x)+ v(x) \big)    \psi_{\h}^2(x )  dx 
=  \lambda_{\h} ,
\label{eq:V2}\end{equation}
where   the eigenfunctions $ \psi_{\h}$ are real and normalized.

\medskip

 {\bf 2.2.}
 As is well known, for homogeneous potentials, the kinetic $K(\psi_{\h})$ and potential $U(\psi_{\h})$ energies are related to the total energy by exact  equalities. 

 \begin{proposition}\label{qclhom}
 Suppose that 
\begin{equation}
   v(tx)  =  t^\alpha v(x),\q t>0,\q \alpha>0.
\label{eq:V3}\end{equation}
 Then
 \begin{equation}
K(\psi_{\h})= \alpha (\alpha+2)^{-1} \lambda_{\h}\q \mathrm{and}
\q U(\psi_{\h})=  2 (\alpha+2)^{-1} \lambda_{\h}.  
\label{eq:V4}\end{equation}
\end{proposition} 

Indeed, it suffices to use that $rv_r(x)=\alpha v(x)$ in \e{eq:V2}.

This result remains true for some $\alpha<0$; however to define $H_{\h}$ as a self-adjoint operator,  we have to require that    $|\alpha|$ be not too large.
 
\medskip

 {\bf 2.3.}
  The case of homogeneous potentials is of course exceptional. In general,
  one cannot expect (if $d>1$) to have even a semiclassical asymptotics   of the
kinetic (or potential) energy. Let us consider a simple

\begin{example}\label{ex}
 Set
 \begin{equation}
v(x)= |x_{1}|^{\alpha_{1}}  + |x_{2}|^{\alpha_{2}},\q x =(x_{1},x_{2})\in{\Bbb R}^2,\q \alpha_{j} >0,\q j=1,2.
\label{eq:vv1}\end{equation}
For every $\lambda>0$ and every
\[
u \in ( 2(\alpha_{1}+2)^{-1}\lambda, 2(\alpha_{2}+2)^{-1}\lambda), 
\]
there exist sequences of eigenvalues $\lambda_{\h}$ and eigenfunctions $\psi_{\h}$ of the operator $H_{\h}$ such that $\lambda_{\h}\to \lambda$ and $U(\psi_{\h})\to u $ as $\h\to 0$.
\end{example}

Indeed, let $ a_{n}^{(j)}$, $n=1,2\ldots$, be eigenvalues of the one-dimensional operators
$-  D_{j}^2+  |x_{j} |^{\alpha_{j}}$,  $j=1,2$.
Eigenvalues $\lambda_{\h}$ of operator \e{eq:I1} with   potential \e{eq:vv1} are given by the formula
\begin{equation}
\lambda_{\h}=\lambda_{\h}^{(1)}+\lambda_{\h}^{(2)}
\label{eq:vv2}\end{equation}
where $\lambda_{\h}^{(j)}=\h^{\gamma_{j} }a_{n_{j}}^{(j)}$,  $j=1,2$,  $\gamma_{j}=2\alpha_{j} (\alpha_{j}+2)^{-1}$ and $n_{j}=1,2\ldots $ are arbitrary.
Normalized eigenfunctions $\psi_{\h}$ are scaled products of normalized  eigenfunctions $\varphi_{n_{j}}^{(j)}$ of these two one-dimensional operators, that is
\[
\psi_{\h}(x)= \h^{-(\beta_{1}+\beta_{2})/2}\varphi_{n_{1}}^{(1)}(\h^{-\beta_{1}}x_{1}) \varphi_{n_{2}}^{(2)}(\h^{-\beta_{2}}x_{2}),\q \beta_{j}=2 (\alpha_{j}+2)^{-1}.
\]
For this eigenfunction, the potential energy equals
\[
U(\psi_{\h}  )=  \h^{\gamma_{1} }\int_{\Bbb R}|x_{1}|^{\alpha_{1}}\varphi_{n_{1}}^{(1)}( x_{1})^2 dx_{1}+   \h^{\gamma_{2} }\int_{\Bbb R}|x_{2}|^{\alpha_{2}}\varphi_{n_{2}}^{(2)}( x_{2})^2 dx_{2}
\]
so that in view of the second formula \e{eq:V4}
\begin{equation}
U(\psi_{\h}  ) 
= \beta_{1}\h^{\gamma_{1} }a_{n_{1}}^{(1)}+  \beta_{2}\h^{\gamma_{2} }a_{n_{2}}^{(2)}.
\label{eq:vv4}\end{equation}

Let us now take into account that
\[
a_{n}^{(j)}= c_{j}n^{\gamma_{j}} (1+o(1)),\q c_{j} > 0,\q n\to\infty.
\]
Pick some numbers $\mu_{j}>0$ and set $n_{j}=[(\mu_{j}c_{j}^{-1})^{1/\gamma_{j}}\h^{-1}]$ where $[b]$ is the integer part of a number $b$. Then 
\begin{equation}
\lambda_{\h}^{(j)}=\h^{\gamma_{j} }a_{n_{j}}^{(j)}\to \mu_{j}
\label{eq:v4x}\end{equation}
as $\h\to 0$. Moreover,  it follows from \e{eq:vv4} that
\begin{equation}
\lim _{\h\ri 0}U(\psi_{\h})=2(\alpha_{1}+2)^{-1}\mu_{1}+ 2(\alpha_{2}+2)^{-1}\mu_{2}.
\label{eq:vv4x}\end{equation}
If $\mu_{1}+\mu_{2}=\lambda$, then according to \e{eq:vv2}  and  \e{eq:v4x}  $\lambda_{\h}\to \lambda$. However the limit \e{eq:vv4x} may take arbitrary values between $2(\alpha_{1}+2)^{-1}\lambda$ and $2(\alpha_{2}+2)^{-1}\lambda$.

\section{Estimates of the kinetic energy}

  {\bf 3.1.}
In addition to the virial theorem, we need  results   on  a decay of eigenfunctions in the classically forbidden region. We suppose that  
\[
\liminf_{|x|\to \infty} v(x) =: v_{\infty}>0 
\]
and consider eigenvalues  $\lambda_{\h}$ which belong to a neighborhood of some non-critical energy $\lambda_{0}\in (0, v_{\infty})$. Our construction works for $v\in C^1 ({\Bbb R}^d)$, but in order to use the  results   on  a decay of eigenfunctions we assume that $v\in C^\infty ({\Bbb R}^d)$.
As always, eigenfunctions $\psi_{\h}$ of the operator $H_{\h}$ are  real and normalized.
Let the sets $F(\lambda )$ and $G(\lambda )$ be defined by the formulas 
\[
F (\lambda )=\{x\in{\Bbb R}^d: v(x) < \lambda \}, \q G(\lambda) ={\Bbb R}^d\setminus F(\lambda).
\]
Then (see \cite{ShHe}, \cite{Simon} as well as \cite{Hel} and references therein) for all {\it  fixed}  $\varepsilon >0$, 
\begin{equation}
 \int _{G ( \lambda_{\h} +\varepsilon )} \big(  \h^2 |\nabla \psi_{\h} |^2 + \psi_{\h}^2(x ) \big)  dx\to 0 
\label{eq:ClFor}\end{equation}
as $\h\to 0$. Actually,    eigenfunctions $\psi_{\h}$ decay exponentially as $\h\to 0$ in the classically forbidden region $G ( \lambda_{\h} +\varepsilon )$, but we do not need this result.

Relation  \e{eq:ClFor} can be supplemented by an estimate of the potential energy. Of course the next lemma is useful only in the case when $v(x)$ is not bounded at infinity.

\begin{lemma}\label{EC}
 For all $\varepsilon>0$, we have
 \[
 \int_{G(\lambda_{\h}+\varepsilon)}
   v(x) \psi_{\h}^2(x)  dx\to 0  
 \]  
 as $  \h\to 0$.   
   \end{lemma} 
   
    \begin{pf} 
    In view of \e{eq:ClFor} it suffices to check that, for some $R$,
   \begin{equation}
\lim_{\h\to 0} \int_{|x|\geq R}
   v(x) \psi_{\h}^2(x)  dx = 0.   
\label{eq:C1}\end{equation}
Choose $R$ such that $v(x)\geq \lambda_{\h}+\varepsilon$ for $|x|\geq R /2$. Let
        $\eta\in C ^\infty ({\Bbb R}^d)$ be such that   $\eta  ( x)\geq 0$,   $\eta  ( x)=0$ for $|x|\leq R/2$ and  $\eta  ( x)=1$ for $|x|\geq R$. Multiplying equation \e{eq:I2}  by $\eta$ and integrating by parts, we see that
    \begin{equation}
   \int_{{\Bbb R}^d}\eta\big( \h^2 |\nabla \psi_{\h} |^2+   v    \psi_{\h}^2 \big) dx=
   \lambda_{\h}   \int_{{\Bbb R}^d}\eta       \psi_{\h}^2   dx
 -  \h^2  \int_{{\Bbb R}^d} \la\nabla \psi_{\h}, \nabla\eta\ra \psi dx .
\label{eq:C2}\end{equation}
Let us consider the right-hand side.
 The first integral  tends to zero as $\h\to 0$ because $\eta(x)=0$ in the classically allowed region.
 Using the Schwarz inequality and  relation \e{eq:Tot}, we estimate the second term by $\h \lambda_{\h}^{1/2} \max |\nabla \eta (x)|$. Therefore expression \e{eq:C2} tends to zero as $\h\to 0$ which proves \e{eq:C1}.
      \end{pf}

 \medskip
 
 {\bf 3.2.}
 Relations \e{eq:V1M} and \e{eq:ClFor} can be combined. The simplest example is given in the next statement.
 
 \begin{proposition}\label{qclexp}
 Let the function $r v_{r}(x) v(x)^{-1}$ be bounded as $|x|\to \infty$. Choose some $\lambda_{0}>0$.
Suppose that, for some $\varepsilon_{0}\in (0,  2 \lambda_{0} )$,        the potential $v(x)$ admits representation \e{eq:V3} in the region $F (\lambda_{0}+\varepsilon_{0})$. Then for $\lambda_{\h}\in  (\lambda_{0}-\varepsilon_{0}/2, \lambda_{0}+\varepsilon_{0}/2 )$,   we have, as $\h\to 0$, the asymptotic relations 
 \begin{equation}
K(\psi_{\h})= \alpha (\alpha+2)^{-1} \lambda_{\h} +o(1) \q \mathrm{and}
\q U(\psi_{\h})=  2 (\alpha+2)^{-1} \lambda_{\h} +o(1) .  
\label{eq:V5}\end{equation}
 \end{proposition}

Indeed, it follows from       \e{eq:V2}  and Lemma~\ref{EC} that
\[
   \int _{F(\lambda_{0}+\varepsilon_{0})} \big( 2^{-1}r v_{r}(x)+ v(x) \big)    \psi_{\h}^2(x ) dx =\lambda_{\h}+o(1) .
\]
By virtue of \e{eq:V3} the integral in the left-hand side equals
 \[
   2^{-1}( \alpha+2) \int _{F(\lambda_{0}+\varepsilon_{0})}   v(x)      \psi_{\h}^2(x )  dx .
   \]
 Therefore using again Lemma~\ref{EC}, we obtain \e{eq:V5}.
 
 Estimates \e{eq:maink} on the kinetic energy or, equivalently, \e{eq:vx} on the potential energy
 can be obtained under much weaker assumptions on $v(x)$.   In view of \e{eq:ClFor} the principal difficulty is to exclude that eigenfunctions $\psi_{\h} (x)$ are localized in   neighborhoods of the surfaces $v(x)=\lambda_{\h}$.   

  \medskip
 
 {\bf 3.3.}
  Let us formulate the main result of this paper.

\begin{theorem}\label{vir1M}
Let $\lambda_{0} > 0$ be a non-critical energy, and let for some $($sufficiently small$)$ $\varepsilon_{0}>0$ and all $\lambda  \in (\lambda_{0}-\varepsilon_{0} , \lambda_{0}+\varepsilon_{0} )$
  \begin{equation}
F (\lambda  ) =\bigcup_{n=1}^N F_{n} (\lambda ), \q F_{n} (\lambda  )\cap F_m (\lambda  )=\emptyset\q {\rm for}\q n\neq m.
\label{eq:W1}\end{equation}
 Suppose that  for   all $n=1,\ldots, N$ and some points $x_{n} \in F_{n} ( \lambda_{0}-\varepsilon_{0} )$  an inequality  
\begin{equation}
  \la x-x_{n}, \nabla v(x)\ra \, \geq  c_{0} >0, \q x\in F_{n} ( \lambda_{0}+\varepsilon_{0})\cap G(\lambda_{0}-\varepsilon_{0} ),    
\label{eq:V1MN}\end{equation}
holds.  
  Denote by $\psi_{\h}$  normalized eigenfunctions $\psi_{\h}$ of the operator $H_{\h}$ corresponding to eigenvalues $\lambda_{\h}$ in a neighborhood $ (\lambda_{0}-\varepsilon_{0}/2, \lambda_{0}+\varepsilon_{0}/2)$ of the point $\lambda_{0}$. Then 
    inequality \e{eq:maink}  is true for sufficiently small $\h$. 
   \end{theorem}
   
The assumptions of this theorem are, actually, very mild. Roughly speaking, 
we suppose that   the classically allowed region $F (\lambda  )$ consists of a finite number of potential wells. Condition \e{eq:V1MN} means that $v(x) $ increases as $x$ passes through the boundary   of $  F_{n}(\lambda  )$. This is consistent with the fact  that $ F_{n}(\lambda  )$ is a potential well of $v(x)$ for the energy $\lambda$.  If $N=1$, then setting $x_{1}=0$, we obtain that inequality \e{eq:V1MN} reduces to the condition $v_{r}(x)\geq  c_{0} >0$.

 We split the proof in a series of simple lemmas.  Let us estimate the potential energy in two different ways. The  first one is quite straightforward.
 
 \begin{lemma}\label{L1}
   For all   $\varepsilon>0$ and all $\d \in (0, \lambda_{\h})$, we have
   \[
\int_{F(\lambda_{\h}+ \varepsilon)}  v (x)\psi_{\h}^2(x)dx\leq \lambda_{\h}+\varepsilon -\d 
\int_{F (\lambda_{\h} -\d) }   \psi_{\h}^2(x)dx .
\]
   \end{lemma}

   \begin{pf} 
  Observe that
\[
\int_{F (\lambda_{\h} -\d)}  v (x)\psi_{\h}^2(x)dx
\leq
 (\lambda_{\h} -\d) \int_{F (\lambda_{\h} -\d)}   \psi_{\h}^2(x)dx
\]
and
\begin{align*}
\int_{ F (  \lambda_{\h}+\varepsilon )
\cap G (\lambda_{\h}-\d  )}  v (x)\psi_{\h}^2(x)dx
&\leq 
(\lambda_{\h}+\varepsilon ) \int_{ G (\lambda_{\h} -\d) }  \psi_{\h}^2(x)dx
\nonumber\\
& =(\lambda_{\h}+\varepsilon ) \big(1- \int_{F (\lambda_{\h} -\d) }   \psi_{\h}^2(x)dx \big).
 \end{align*}
   So it suffices to put  these two   estimates together.
 \end{pf}
 
 Combining Lemmas~\ref{EC} and \ref{L1}, we see that for all $\varepsilon>0$ and all  $\d \in (0, \lambda_{\h}) $,
    \begin{equation}
    U (\psi_{\h}) \leq \lambda_{\h}+\varepsilon -\d
\int_{F (\lambda_{\h} - \d) }   \psi_{\h}^2(x)dx + \sigma( \varepsilon , \h) 
\label{eq:VV1}\end{equation}
where   $\sigma( \varepsilon , \h) \to 0$ as $\h\to 0$ if   $\varepsilon$ is {\it fixed}. The notation $\sigma( \varepsilon , \h)  $ will also be used below.

 The second estimate relies on the virial theorem which we need in the following form.

 \begin{lemma}\label{L2}
Let $x_{n}\in F_{n} (\lambda_{0}+ \varepsilon_{0})$ be arbitrary points. Then
\begin{equation}
\sum_{n=1}^N \int _{F_{n} (\lambda_{0}+ \varepsilon_{0})} \Bigl( 2 \h^2 | \nabla  \psi_{\h} (x)|^2  -
 \la x-x_{n}, \nabla v(x) \ra \psi_{\h}^2(x)   \Bigr)dx= o(1), \q \h \to 0.
\label{eq:V1Mcl1}\end{equation}  
   \end{lemma}

 \begin{pf}
Let us use Theorem~\ref{GVT} for a suitable function $a(x)$ which we construct now. 
 Choose functions $\varphi_{n}\in C_{0}^\ii ({\Bbb R}^d)$ such that $\varphi_{n}(x)=1$  for  $x\in F_{n} (\lambda_{0}+ \varepsilon_{0})$ and $\varphi_{n}(x)=0$   away from some neighborhoods of $  F_{n} (\lambda_{0}+ \varepsilon_{0})$ so that  $\supp \varphi_{n}\cap\supp \varphi_m =\emptyset$ if $n\neq m$. We   define the function $a(x)$ by the equality
\begin{equation}
a(x)=  \sum_{n=1}^N | x- x_n|^2 \varphi_n(x).
\label{eq:V1MN1}\end{equation}

Neglecting in \e{eq:V1M} the classically forbidden region, we see that
\begin{equation}
  \int _{F  (\lambda_{0}+ \varepsilon_{0})} \Bigl( 4 \h^2\sum_{j,k=1}^d    a_{jk} \partial_{k}\psi_{\h} \partial_{j}\psi_{\h}   -
 \h^2  (\Delta^2 a) \psi_{\h}^2  -
2  \la\nabla a, \nabla v\ra \psi_{\h}^2   \Bigl)dx= o(1)
\label{eq:V1Mcl}\end{equation}
as $\h\to 0$. 
If $x\in F_{n} (\lambda_{0}+ \varepsilon_{0})$, then   according to \e{eq:V1MN1} we have $(\nabla a)(x)= 2(x - x_{n})$, $a_{jj}(x)=2$ and $a_{jk}(x)=0$ if $j\neq k$.
 Thus, relation \e{eq:V1Mcl1} follows from \e{eq:V1Mcl}.
 \end{pf}

 \begin{lemma}\label{L3}
 Let assumption \e{eq:V1MN} hold for some points $x_{n} \in F_{n} ( \lambda_{0}-\varepsilon_{0} )$, and  set
\begin{equation}
 c_{1}=  \max_{n}\sup_{x\in F_{n} (\lambda_{0} - \varepsilon_{0})} (-\la x-x_{n}, \nabla v(x)\ra ).
\label{eq:c1}\end{equation}
Then
\begin{equation}
2  \h^2 \int_{ {\Bbb R}^d}  | \nabla  \psi_{\h} (x)|^2dx\geq  c_{0} - c_{2}\int_{F(\lambda_0-\varepsilon_{0})} \psi_{\h}^2(x)dx + \sigma (\varepsilon_{0}, \h)
\label{eq:vir8}\end{equation}
where $c_2=c_{0}+c_{1}$ and $\sigma (\varepsilon_{0}, \h)\to 0$
as $\h \to 0$.     
   \end{lemma}
   
    \begin{pf}
   According to \e{eq:V1MN} and \e{eq:c1}, we have
 \begin{align*}
    \int _{F_{n} (\lambda_{0}+ \varepsilon_{0})}    \la x- x_{n}, \nabla v (x)\ra \psi_{\h}^2 (x)  dx   \geq &  c_{0}
        \int _{F_{n} (\lambda_{0}+ \varepsilon_{0}) \cap G(\lambda_0- \varepsilon_{0} )}    \psi_{\h}^2 (x)  dx
        \\
        & -   c_{1}
           \int _{  F_{n}(\lambda_0- \varepsilon_{0} )}    \psi_{\h}^2 (x)  dx .
    \end{align*}
         Summing these estimates over $n=1,\cdots, N$ and using \e{eq:V1Mcl1}, we see that
           \begin{align*}
2   \h^2 \int_{ F(\lambda_0+ \varepsilon_{0} ) }  | \nabla  \psi_{\h} (x)|^2dx
  \geq & c_0 \int_{F  (\lambda_{0}+ \varepsilon_{0}) \cap G(\lambda_0- \varepsilon_{0} )}\psi_{\h}^2(x)dx 
  \\
  &- c_{1} \int_{F(\lambda_0-\varepsilon_{0})} \psi_{\h}^2(x)dx + o(1)
       \end{align*}
as $\h\to 0$.
The first integral in the right-hand side equals $1$ minus the integrals of $\psi_{\h}^2(x)$ over $F(\lambda_0-\varepsilon_{0})$ and $G(\lambda_0+\varepsilon_{0})$. The   integral over $G(\lambda_0+\varepsilon_{0})$ tends to zero because $ G(\lambda_0+\varepsilon_{0})$ lies in the classically forbidden region. 
           \end{pf}

      Using         the energy conservation  \e{eq:Tot} and the obvious inclusion $F(\lambda_0-\varepsilon_{0})\subset F(\lambda_{\h}-\varepsilon_{0}/2)$, we deduce from     \e{eq:vir8} the estimate
       \begin{equation}
 U (\psi_{\h}) \leq  \lambda_{\h} - 2^{-1} c_{0} + 2^{-1}  c_{2}\int_{F(\lambda_{\h}-\varepsilon_{0}/2)} \psi_{\h}^2(x)dx + \sigma (\varepsilon_{0}, \h).
\label{eq:VV2}\end{equation}
  If $c_{2}\leq 0$, then \e{eq:VV2} directly implies  \e{eq:vx}. So below we assume $c_{2}> 0$.

    Now we are in a position to conclude the {\it proof} of Theorem~\ref{vir1M}. Let us   compare estimates \e{eq:VV1} where we set $\d=\varepsilon_{0}/2$ and \e{eq:VV2}. Roughly speaking, if the integral
\[
X=\int_{F(\lambda_{\h}-\varepsilon_{0}/2)}\psi_{\h}^2(x) dx
\]
  is small, then we use \e{eq:VV2}.  If it is big, we use \e{eq:VV1}. To be more precise, estimates \e{eq:VV1} and \e{eq:VV2} imply that
  \[
U (\psi_{\h}) 
\leq \lambda_{\h}+ 2^{-1} \max_{0\leq X\leq 1}\min\{ 2 \varepsilon -  \varepsilon_{0} X,-c_{0}+c_2 X\}+ \sigma (\varepsilon, \varepsilon_{0}, \h),
\]
where $\sigma (\varepsilon, \varepsilon_{0}, \h)\to 0$ as $\h\to 0$
for    fixed $\varepsilon $ and  $\varepsilon_{0} $.
Observe   that
 \[
\max_{0\leq X\leq 1}\min\{ 2 \varepsilon -  \varepsilon_{0} X,-c_{0}+c_2 X\}\leq 2 \varepsilon -  \varepsilon_{0} c_{0}(1+c_{2})^{-1}
\]
(if $\varepsilon_{0}\leq 1$). Since $\varepsilon$ is arbitrary   small, this yields estimates \e{eq:vx} and hence \e{eq:maink}. In these estimates $c$ is any number smaller than
$2^{-1}\varepsilon_{0} c_{0}(1+c_0+ c_{1})^{-1}$.
 
  \medskip
 
 {\bf 3.4.}
  Our lower bound on the potential energy (and hence an upper bound on the kinetic energy) is almost trivial.
 
   \begin{proposition}\label{vir1P}
 Let,   for some $ \lambda_{0}>0$ and $\varepsilon_{0}>0$, representation \e{eq:W1} hold for $\lambda=\lambda_{0} + \varepsilon_{0}$.
 Suppose that for all $n=1,\ldots, N$ there exist points $x_{n} \in F_{n} (\lambda_{0})$ such that the estimates
 \begin{equation}
  \la x-x_{n}, \nabla v(x)\ra \, \leq c_{0} v(x), \q x\in F_{n} (\lambda_{0}),     
\label{eq:V8x}\end{equation}
are satisfied with some constant $c_{0}>0$. Then
  \begin{equation}
U (  \psi_{\h})  \geq c \lambda_{\h}   
\label{eq:V4c}\end{equation}
 for all $\lambda_{\h} \in (\lambda_{0}- \varepsilon_{0}/2, \lambda_{0} +\varepsilon_{0}/2)$, an arbitrary $c < 2 ( c_{0}+ 2 )^{-1}$ and sufficiently small $ \h$.
   \end{proposition}
   
   \begin{pf} 
Comparing  relation  \e{eq:V1Mcl1} with assumption \e{eq:V8x}, we see that
      \[
  2 \h^2  \int _{F(\lambda_{0}+ \varepsilon_{0})}    | \nabla  \psi_{\h} (x)|^2  dx \leq
  c_{0} \int _{F(\lambda_{0}+ \varepsilon_{0})}  v(x)    \psi_{\h}^2 (x)   dx +\sigma (\varepsilon_{0},\h).
\]
Then using relation   \e{eq:ClFor} and Lemma~\ref{EC}, we obtain the  estimate
\[
 2 K(\psi_{\h}) \leq c_{0}  \, U(\psi_{\h} ) + o (\h).
\]
In view of the energy conservation \e{eq:Tot}, this yields \e{eq:V4c}.
 \end{pf}

 Assumption \e{eq:V8x} 
  essentially means that, inside every well $F_{n} (\lambda_{0})$,   the function $v(x)$ may equal zero only at the point $x_{n}$. If, for example, 
  \[
   v(x)= v_{n}|x-x_{n}|^{\alpha_{n}}, \q v_{n}>0, \q \alpha_{n}>0, \q x\in F_{n} (\lambda_{0}), \q n=1,\ldots, N,
   \]
   then estimate \e{eq:V8x} holds with $c_{0}=\max\{\alpha_{1},\ldots, \alpha_{N}\}$.

 \medskip
 
 I thank B. Helffer and D. Robert for   useful discussions.


      \end{document}